\documentclass[12pt,a4paper]{article}
\usepackage{amsfonts, amsmath, amssymb}

\newtheorem{prop}{Proposition}
\newtheorem{lemma}{Lemma}
\newtheorem{theorem}{Theorem}
\newtheorem{rem}{Remark}
\newtheorem{exmp}{Example}
\newtheorem{cor}{Corollary}

\author{Mark Pankov}
\title{Transformations of Grassmannians and automorphisms of
classical groups}
\date{}
\begin{document}

\maketitle

\begin{abstract}
\parskip=3mm
\parindent=0mm

\noindent We consider transformations preserving certain linear
structure in Grassmannians and give some generalization of the
Fundamental Theorem of Projective Geometry and the Chow Theorem
[Ch]. It will be exploited to study linear $(k,n-k)$-involutions,
$1<k<n-1$. The analogy of the J. Dieudonn\'{e} and C. E. Rickart
result will be obtained.


{\bf MSC-class}: 14M15, 14L35, 20G15, 51N30.
\end{abstract}

\section{Introduction}

Let $V$ be an $n$-dimensional vector space over
some field $F$. Denote by ${\mathbb G}_{k}(V)$
the Grassmannian consisted of $k$-dimensional linear
subspaces of $V$. In what follows a
$k$-dimensional linear subspace
will be called a $k$-dimensional plane if
$k>1$ and a line if $k=1$.

\subsection{The Fundamental Theorem of Projective
Geometry and Chow Theorem}

Let us consider a transformation $f$
(bijection onto itself) of the Projective space
${\mathbb G}_{1}(V)$ such that
$f$ and the inverse transformation $f^{-1}$
map each collection of linearly independent lines
to a collection of linearly independent lines
(we say that lines are linearly independent
if non-zero vectors lying on them are
linearly independent).
Each collineation (semilinear transformation) of $V$
induces the transformation of ${\mathbb G}_{1}(V)$
satisfying this condition
(note that two transformations of
${\mathbb G}_{k}(V)$ induced by collineations
$f$ and $f'$ are coincident if and only if
$f'=af$ for some $a\in F$).

The inverse statement is known as
the Fundamental Theorem of Projective Geometry.
It can be formulated in the following form: if $n\ge 3$ then each
transformation of ${\mathbb G}_{1}(V)$ satisfying the condition
considered above is induced by a collineation  (see, for example,
[D2, O'M1]).  It must be pointed out that this result was
successfully exploited to description of automorphisms of linear
groups [D1, D2, O'M1, O'M2, R].

An analogy of the Fundamental Theorem of
Projective Geometry for  transformations
of ${\mathbb G}_{k}(V)$, $k>1$,
was given by W. L. Chow [Ch] (see also [D2]).
The formulation of the Chow Theorem is based on the notion
of the distance between planes.

For two arbitrary $k$-dimensional planes $S,S'\subset V$
the number
$$k-\dim S\cap S'$$
is called the {\it distance} between $S$ and $S'$.
It is easy to see that the distance
is equal to the smallest number $i$ such that
$S$ and $S'$ are contained in some
$(k+i)$-dimensional plane.
Two planes
will be called {\it adjacent} if the distance between
them is equal to $1$.

We consider transformations of Grassmannians
preserving the distance between planes.
A transformation $f$ of ${\mathbb G}_{k}(V)$ preserves the
distance if and only if $f$ and $f^{-1}$ map
any two adjacent planes to adjacent planes
(it is trivial, see [D2]). If
$k=1$ or $n-1$ then any two elements of ${\mathbb G}_{k}(V)$ are
adjacent and each transformation of ${\mathbb G}_{k}(V)$
preserves the distance.

Clearly, each
transformation of ${\mathbb G}_{k}(V)$
induced by a collineation
preserves the distance.
The Chow Theorem states that
if $n\ge 3$, $1<k<n-1$ and $n\ne 2k$ then the inverse statement
holds true: a transformation of ${\mathbb G}_{k}(V)$
preserving the distance between planes is induced by some
collineation.
For the case when $n=2k$ it fails.

Consider a non-degenerate sesquilinear form $\Omega$ on $V$.
It defines the bijection $f_{k\,n-k}(\Omega)$
of the Grassmannian ${\mathbb G}_{k}(V)$ onto
the Grassmannian ${\mathbb G}_{n-k}(V)$
which transfers each plane to the
$\Omega$-orthogonal complement.
This bijection
preserves the distance.
The bijections defined by forms
$\Omega$ and $\Omega'$ are coincident
if and only if
there exists $a\in F$ such that $\Omega'=a\Omega$.

If $3\le n=2k$ then each $f_{k\,k}(\Omega)$ is a transformation of
${\mathbb G}_{k}(V)$ which is not induced by a collineation.
For this case the Chow Theorem states that
each transformation of ${\mathbb G}_{k}(V)$
preserving the distance between planes is induced by a
collineation or defined by
a non-degenerate sesquilinear form.

The similar statements were also proved for some
homogeneous spaces, for example, for the space of
all null planes of a non-degenerate symplectic form
(see [Ch]).
Other generalization of the
Fundamental Theorem of Projective Geometry and
the Chow Theorem was proposed by  J. Tits [T].

\subsection{Regular transformations of Grassmannians}

Now we introduce so-called $R$-subsets of
Grassmannians. They will be used to formulate a
statement generalizing the results considered
above.

We say that ${\cal R} \subset {\mathbb G}_{k}(V)$ is
an $R$-{\it set} if
there exists a base for $V$ such that
each plane belonging to ${\cal R}$ contains
$k$ vectors from this base; in other words,
elements of ${\cal R}$
are coordinate planes for some coordinate system
for $V$.  Any base and any coordinate system satisfying
that condition will be called {\it associated}
with ${\cal R}$.

A coordinate system for $V$ has
$$
\binom{n}{k}
:=\frac{n!}{k!(n-k)!}
$$
distinct $k$-dimensional coordinate planes.
Therefore an $R$-subset of ${\mathbb G}_{k}(V)$ contains at most
$\binom{n}{k}$ elements.
An $R$-set ${\cal R}$ will be called  {\it
maximal} if any $R$-set containing
${\cal R}$ coincides with it.
Each
$R$-set is contained in some maximal $R$-set and
an $R$-subset of ${\mathbb G}_{k}(V)$ is maximal
if and only if it contains $\binom{n}{k}$ elements.
Thus a maximal $R$-subset of ${\mathbb G}_{k}(V)$
consists of all $k$-dimensional coordinate
planes for some coordinate system.

An $R$-subset of ${\mathbb G}_{1}(V)$ is
a collection of linearly independent lines;
$R$-subsets of ${\mathbb G}_{n-1}(V)$ are known
as {\it arrangements} [O].
For the general case $R$-subsets of Grassmannians are
discrete sets which
can be considered as a generalization of
collections of linearly independent lines.

A transformation $f$ of ${\mathbb G}_{k}(V)$
will be called  {\it regular} if
$f$ and $f^{-1}$
preserve the class of $R$-sets.
Any transformation induced by
a collineation is regular.

We say that a bijection $f$ of ${\mathbb G}_{k}(V)$ onto
${\mathbb G}_{n-k}(V)$ is {\it regular} if $f$ and $f^{-1}$
transfer
each $R$-set to
an $R$-set.
It is easy to see that
for any non-degenerate sesquilinear form $\Omega$ on $V$
the bijection $f_{k\,n-k}(\Omega)$ is regular.
Note that for the case when $m\ne k,n-k$
there are not regular bijections of ${\mathbb G}_{k}(V)$ onto
${\mathbb G}_{m}(V)$ since
the equality
$$\binom{n}{k}=\binom{n}{m}$$
holds only for $m=k,n-k$.

If $n=2$ then
any two lines generate a maximal
$R$-subset of ${\mathbb G}_{1}(V)$
and each transformation of ${\mathbb G}_{1}(V)$ is
regular. For the case when $n\ge 3$ the class
of regular transformations of ${\mathbb G}_{k}(V)$
consists only of the transformations considered above.
\begin{theorem}
If $n\ge 3$ then the following two statements hold true:
\begin{enumerate}
\item[---] for the case when $n\ne 2k$ all regular
transformations of ${\mathbb G}_{k}(V)$ are induced
by collineations ;
\item[---] if $n=2k$ then a regular
transformation of ${\mathbb G}_{k}(V)$ is induced
by a collineation or
defined by some non-degenerate sesquilinear form.
\end{enumerate}
\end{theorem}
\begin{cor}
If $n\ge 3$ and $n\ne 2k$ then
each regular bijection
of ${\mathbb G}_{k}(V)$ onto
${\mathbb G}_{n-k}(V)$
is defined by
a non-degenerate sesquilinear form.
\end{cor}
For $k=n-1$ Theorem 1.1 is a simple
consequence of
the Fundamental Theorem of Projective Geometry.
For the general case it will be proved in the next section.
Our proof will be based on properties of
sertain number characteristic of $R$-sets, so-called, the degree
of inexactness (Subsections 2.1 and 2.2).

\subsection{Linear involutions and automotphisms of classical
groups}

Let $\sigma\in {\mathfrak G}{\mathfrak L}(V)$
be an involution (i.e. $\sigma^{2}=Id$) and the characteristic of
the field $F$ is not equal to $2$. Then there exist
two
invariant planes $U_{+}(\sigma)$ and $U_{-}(\sigma)$
such that
$$\sigma(x)=x\;\mbox{ if }\;x\in U_{+}(\sigma)\,,
\;\;\;\sigma(x)=-x\;\mbox{ if }\;x\in U_{-}(\sigma)$$
and
$$V=U_{+}(\sigma)+U_{-}(\sigma)\,.$$
We say that
$\sigma$ is a $(k,n-k)$-{\it involution} if the dimensions
of the planes $U_{+}(\sigma)$ and $U_{-}(\sigma)$
are equal to $k$ and $n-k$, respectively.

The set of all $(k,n-k)$-involutions will
be denoted by ${\mathfrak I}_{k\,n-k}(V)$.  If the number $n-k$
is even then $(k,n-k)$-involutions generate the special linear
group ${\mathfrak S}{\mathfrak L}(V)$.  For the case when $n-k$
is odd they generate the group of all  linear transformations
$f\in {\mathfrak G}{\mathfrak L}(V)$ satisfying the condition
$\det f=\pm 1$.

In what follows we will consider transformations
of ${\mathfrak I}_{k\,n-k}(V)$ preserving
the commutativity. For example, for each
collineation $g$ of $V$ the transformation
$$\sigma \to g\sigma g^{-1}$$
of ${\mathfrak I}_{k\,n-k}(V)$ satisfies
this condition. Now consider
a correlation $h$
(a semilinear bijection of $V$ onto the dual space $V^{*}$).
Then the transformation
$$
\sigma \to h^{-1}\check{\sigma}h
$$
of ${\mathfrak I}_{k\,n-k}(V)$
(here $\check{\sigma}$ is the contragradient)
preserves
the commutativity.

It was proved by
J. Dieudonn\'{e} [D1] and C. E. Rickart [R]
that {\it for the cases $k=1,n-1$ any transformation of
${\mathfrak I}_{k\,n-k}(V)$ preserving the commutativity
is defined by a collineation or a correlation}.
The proof was based on
the Fundamental Theorem of Projective Geometry and
G. W. Mackey's results [M].
It is not difficult to see that automorphisms
of the groups ${\mathfrak S}{\mathfrak L}(V)$ and
${\mathfrak G}{\mathfrak L}(V)$
induce transformations of the set
${\mathfrak I}_{k\,n-k}(V)$
preserving the commutativity.
Therefore the well-known
description of automorphisms of the groups
${\mathfrak G}{\mathfrak L}(V)$ and
${\mathfrak S}{\mathfrak L}(V)$ is a simple
consequence of the theorem given
above (see [D1, D2, R]).

In Section 3
the similar statement will be proved
for transformations of
${\mathfrak I}_{k\,n-k}(V)$, $1<k<n-1$.

\section{Proof of Theorem 1.1}

\setcounter{theorem}{0}
\setcounter{equation}{0}

\subsection{Degree of inexactness}

It was noted above that each $R$-set
is contained in some maximal $R$-set.
We say that an $R$-set is {\it exact} if
there exists unique maximal $R$-set containing it.
Certainly any maximal $R$-set is exact.

For an $R$-set ${\cal R}$
consider an exact $R$-set
${\cal R}'$
containing ${\cal R}$ and
having the minimal number of elements (i.e.
such that the inequality
$|{\cal R}'|\le |{\cal R}''|$ holds
for any other exact $R$-set
${\cal R}''$ containing ${\cal R}$).
The number
$$ \deg({\cal R}):=|{\cal R}'|-|{\cal R}|$$
will be called the {\it degree of inexactness} of
${\cal R}$.
An $R$-set is exact if and only if the degree of inexactness
is zero.

It is trivial that the following statement
holds true.
\begin{lemma}
Regular transformations of
${\mathbb G}_{k}(V)$
and regular bijections of ${\mathbb G}_{k}(V)$ onto
${\mathbb G}_{n-k}(V)$
preserve the degree of inexactness.
\end{lemma}
If $k=1$ then for any $R$-set
${\cal R}\subset{\mathbb G}_{k}(V)$ we have
$\deg({\cal R})= n-|R|$ and our
$R$-set is exact if and only if
it is maximal.
Lemma 2.1 guarantees the fulfilment of the
similar statement for $k=n-1$.
For the general case there exist exact
$R$-sets which are not maximal.

Now we give a few technical definitions
which will be used in what follows.

For an $R$-set ${\cal R}'\subset {\mathbb G}_{k}(V)$
fix an exact $R$-set
${\cal R}''$ containing ${\cal R}'$
and satisfying the condition
\begin{equation}
\deg({\cal R}')=|{\cal R}''|-|{\cal R}'|\,.
\end{equation}
There exists  unique maximal $R$-set
${\cal R}$ containing ${\cal R}''$.
Let ${\cal R}^{m}$ be the set of all
$m$-dimensional coordinate planes for the coordinate system
associated with ${\cal R}$ (this coordinate system is uniquely
defined since ${\cal R}$ is a maximal $R$-set).
The set ${\cal R}^{m}$ will be called the maximal
$R$-subset of ${\mathbb G}_{m}(V)$ {\it associated} with ${\cal
R}$.  For a plane $S$ belonging to ${\cal
R}^{m}$ denote by ${\cal R}(S)$ the set of all planes $U\in {\cal
R}$ incident to $S$, i.e. such that
$$U\subset S\; \mbox{ if }\; k<m$$
and
$$S\subset U\; \mbox{ if }\; k>m\,.$$
The set ${\cal R}(S)$ contains $\binom{m}{k}$
elements if $m>k$ and $\binom{n-m}{k-m}$ elements if $m<k$.

Let $L_{1},...,L_{n}$ be lines generating
the set ${\cal R}^{1}$.
For each $i=1,...,n$
$${\cal R}'_{i}:={\cal R}'\cap {\cal R}(L_{i})$$
is the  set of all planes belonging to ${\cal R}'$
and containing the line $L_{i}$.
Let
$S_{i}$ be the intersection of all
planes belonging to ${\cal R}'_{i}$.
Then
$${\cal R}'_{i}={\cal R}'\cap {\cal R}(S_{i})\,.$$
The dimension of the plane $S_{i}$ will be denoted
by $n_{i}$; for the case when the set
${\cal R}'_{i}$ is empty we write $n_{i}=0$.
Denote by $n({\cal R}')$ the
number of all $i$ such that $n_{i}=1$.
It is not difficult to see that
\begin{enumerate}
\item[---]{\it the number $n({\cal R}')$ does not depend on the
choice of an exact $R$-set ${\cal R}''$
containing ${\cal R}'$ and satisfying the
condition} (2.1),
\item[---] {\it the $R$-set ${\cal R}'$ is exact if and
only if $n({\cal R}')=n$}.
\end{enumerate}
Let us consider two examples.
\begin{exmp}{\rm
Suppose that the $R$-set ${\cal R}'$
coincides with ${\cal R}(L_{j})$ and
$k\ge n-k$. Then $n_{j}=1$.
If $i\ne j$ then $S_{i}$ is the two-dimensional
plane containing the lines $L_{j}$ and
$L_{i}$; therefore $n_{i}=2$ and $n({\cal R}')=1$.
Consider a set
$$\{i_{1},...,i_{k}\}\subset\{1,...,n\}\setminus\{j\}$$
and the plane $U'\in{\cal R}$
containing the lines $L_{i_{1}}$,...,$L_{i_{k}}$.
For each $p=1,...,k$ the intersection of $U'$
with $S_{i_{p}}$ coincides with the line $L_{i_{p}}$.
This implies the equality
$$n({\cal R}'\cup\{U'\})=k+1\,.$$
If $k=n-1$ then the $R$-set
${\cal R}'\cup\{U'\}$ is exact and
$\deg{\cal R}'=1$.
For the case when $k<n-1$
the set
$$\{1,...,n\}\setminus\{j,i_{1},...,i_{k}\}$$
contains less than $k$ elements.
Denote them by $j_{1},...,j_{q}$ and consider
the plane $U''\in {\cal R}$
containing the lines $L_{j_{1}}$,...,$L_{j_{q}}$.
If $U''$ does not contain $L_{j}$ then
the intersections of $U''$ with the planes
$S_{j_{1}}$,...,$S_{j_{q}}$ are the lines
$L_{j_{1}}$,...,$L_{j_{q}}$.
Thus the $R$-set ${\cal R}'\cup\{U',U''\}$
is exact and $\deg{\cal R}'=2$.
}\end{exmp}
\begin{exmp}{\rm
Now suppose that $k\le n-k$ and the $R$-set ${\cal R}'$
coincides with ${\cal R}(S)$, where $S$
is a plane belonging to
${\cal R}^{n-1}$.
Consider unique line $L_{i}$ which is not
contained in $S$. Then $n_{i}=0$
and $n_{j}=1$ for each $j\ne i$, thus
$n({\cal R}')=n-1$.
For a non-degenerate sesquilinear form
$\Omega$ the bijection $f_{k\,n-k}(\Omega)$ transfers
${\cal R}'$ to the set considered in Example 2.1.
Then Lemma 2.1 guarantees the fulfilment of the equality
$\deg{\cal R}'=2$ for the case when $k>1$.
}\end{exmp}
\begin{theorem}
For the case when $1<k<n-1$
the following three statements hold true:
\begin{enumerate}
\item [{\rm(i)}]
if $n-k<k$ and the set ${\cal R}'$ contains not less than
$\binom{n-1}{k-1}$
elements then $\deg({\cal R}') \le 2$
and the equality $\deg({\cal R}')=2$
holds if and only
${\cal R}'$ is the set considered in Example {\rm 2.1};
\item [{\rm(ii)}]
if $k<n-k$ and the set ${\cal R}'$ contains not
less than
$\binom{n-1}{k}$
elements then $\deg({\cal R}') \le 2$ and the
equality $\deg({\cal R}')=2$ holds if and only if
${\cal R}'$ is the set considered in Example {\rm 2.2};
\item [{\rm(iii)}]
if $n=2k$ and the set ${\cal R}'$ contains not less than
$\binom{n-1}{k}=\binom{n-1}{k-1}$
elements then $\deg({\cal R}') \le 2$ and the equality $\deg({\cal R}')=2$
holds if and only
${\cal R}'$ is one of the sets considered in Examples {\rm 2.1} and
{\rm 2.2}.
\end{enumerate}
\end{theorem}

\subsection{Proof of Theorem 2.1}

We start with a few lemmas.
In the second part of the subsection they will
be exploited to prove Theorem 2.1.
Lemma 2.1 shows that
we can restrict
ourself only to the case when
$n-k \le k <n-1$ and the
$R$-set ${\cal R}'$ contains not less than
$\binom{n-1}{k-1}$ elements.
\begin{lemma}
If the condition $n_{i}=0$ holds for some number $i$
then $n=2k$ and
the set ${\cal R}'$ coincides with ${\cal R}(S)$,
where $S$ is a
plane belonging to ${\cal R}^{n-1}$.
\end{lemma}
{\bf Proof.}
Consider the plane $S\in {\cal R}^{n-1}$
which does not contain the line $L_{i}$.
The condition $n_{i}=0$ shows
that the set ${\cal R}'_{i}$ is empty.
Then we have the inclusion
${\cal R}' \subset {\cal R}(S)$
showing that
$$\binom{n-1}{k-1}
\le |{\cal R}'| \le |{\cal R}(S)|=
\binom{n-1}{k}\,.$$
However $\binom{n-1}{k-1}\ge \binom{n-1}{k}$
if $k\ge n-k$
and the inequality can be replaced by
an equality if and only if $n=2k$.  We get the required.
$\square$
\begin{lemma}
The inequality $n_{i} \le n-k$ holds  for
each $i=1,...,n$.
\end{lemma}
{\bf Proof.}
The case $n_{i}=0$ is trivial.
For the case when $n_{i}>0$ the set ${\cal R}'_{i}$ is not empty
and there exists a plane $U\in {\cal R}'$ containing the line
$L_{i}$.  Then $n_{i} \le k$ and
the required inequality holds for $n=2k$.

Let $n-k<k$. Denote by $S$ the
plane belonging to ${\cal R}^{n-1}$ and such that
the line $L_{i}$ is not contained in $S$.
For each plane $U\in {\cal R}'$
the following two cases can be realized:
\begin{enumerate}
\item[---] $U$ contains $L_{i}$ then
$U\in {\cal R}'_{i}\subset {\cal R}(S_{i})$;
\item[---] $U$ does not contain $L_{i}$
then $U\in{\cal R}(S)$.
\end{enumerate}
This implies the inclusion
${\cal R}'\subset{\cal R}(S)\cup{\cal R}(S_{i})$.
Thus
\begin{equation}
\binom{n-1}{k-1}
\le |{\cal R}'|
\le |{\cal R}(S)| + |{\cal R}(S_{i})|=
\binom{n-1}{k}+\binom{n-n_{i}}{k-n_{i}}\,.
\end{equation}
If $n_{i} \ge n-k+1$ then
$$\binom{n-n_{i}}{k-n_{i}}
\le \binom{k-1}{2k-n-1}$$ (see Remark 2.1, it will be given after
the proof) and the inequality (2.2) can be rewritten in the
following form
\begin{equation}
\binom{n-1}{k-1}-\binom{n-1}{k}
\le
\binom{k-1}{2k-n-1}\,.
\end{equation}
We have
\begin{equation*}
\begin{split}
\binom{n-1}{k-1}-\binom{n-1}{k}=\;&
\frac{(n-1)!}{(k-1)!(n-k)!}-\frac{(n-1)!}{k!(n-k-1)!}=
\frac{(n-1)!(2k-n)}{k!(n-k)!}\\
=\;&(2k-n)\underbrace{k(k+1)...(n-2)(n-1)}_{n-k}
\frac{(k-1)!}{k!(n-k)!}
\end{split}
\end{equation*}
and
\begin{equation*}
\begin{split}
\binom{k-1}{2k-n-1}=\;&\frac{(k-1)!}{(2k-n-1)!(n-k)!}=\\
&(2k-n)\underbrace{(2k-n+1)...(k-1)k}_{n-k}
\frac{(k-1)!}{k!(n-k)!}\,.
\end{split}
\end{equation*}
The condition $n-k < k<n-1$ shows that
$$(2k-n)\underbrace{k(k+1)...(n-2)(n-1)}_{n-k} >
(2k-n)\underbrace{(2k-n+1)...(k-1)k}_{n-k}$$
and the inequality (2.3) does not hold.
The fulfilment of the inequality
$n_{i}\le n-k$ is proved.
$\square$
\begin{rem}{\rm
An immediate verification
shows that
$$
\binom{n-k_{1}}{k-k_{1}}
\ge
\binom{n-k_{2}}{k-k_{2}}
$$
for any two natural numbers $k_{1}$ and $k_{2}$
satisfying the condition $0\le k_{1}\le k_{2}\le k$.
}\end{rem}
\begin{lemma}
If there exists a number $i$ satisfying the
condition $n_{i} \ge 3$ then $n_{j}=1$ for each
$j \ne i$.
\end{lemma}
{\bf Proof.}
Consider the
planes $S\in {\cal R}^{n-1}$ and
$S'\in {\cal R}^{n-2}$ such that
$S$ does not contain the line $L_{i}$ and
$S'$ does not contain the lines $L_{i}$ and
$L_{j}$. Then for each plane $U$
belonging to ${\cal R}'$ the following three cases can be
realized.
\begin{enumerate}
\item[---] $U\in {\cal R}'_{i}$.
\item[---] If $U \notin {\cal R}'_{i}$ and $U \in {\cal R}'_{j}$
then $U$ does not contain the line $L_{i}$ and
$U \in {\cal R}(S)$.  The condition $U \in {\cal R}'_{j}$ shows
that $U$ belongs to the set
${\cal R}(S)\cap {\cal R}(S_{j})$.
\item[---] If $U\notin {\cal R}'_{i}\cup{\cal R}'_{j}$ then $U$
does not contain the lines $L_{i}$ and $L_{j}$.
For this case $U\in {\cal R}(S')$.
\end{enumerate}
In other words, we have the inclusion
$$ {\cal R}'\subset {\cal R}(S_{i})\cup
({\cal R}(S)\cap {\cal R}(S_{j}))\cup {\cal R}(S')$$
showing that
\begin{equation}
|{\cal R}'| \le
|{\cal R}(S_{i})|+ |{\cal R}(S)\cap {\cal R}(S_{j})|+
|{\cal R}(S')|\,.
\end{equation}
The set
${\cal R}(S)\cap {\cal R}(S_{j})$
is not empty if and
only if the plane $S_{j}$ is contained in $S$.
It is not difficult to see that
in this case our set contains
$\binom{n-n_{j}-1}{k-n_{j}}$
elements.
Then by (2.4)
$$\binom{n-1}{k-1}\le
\binom{n-n_{i}}{k-n_{i}}+\binom{n-n_{j}-1}{k-n_{j}}+\binom{n-2}{k}\,.$$
For the case when
\begin{equation}
n_{i} \ge 3\;\mbox{ and }\;n_{j} \ge 2
\end{equation}
we have
$$\binom{n-n_{i}}{k-n_{i}}\le \binom{n-3}{k-3}
\;\mbox{ and }\;
\binom{n-n_{j}-1}{k-n_{j}}\le \binom{n-3}{k-2}$$
(Remark 2.1); then
$$
\binom{n-1}{k-1} \le \binom{n-3}{k-3}+ \binom{n-3}{k-2} +
\binom{n-2}{k}\,.
$$
We use the equality
$$\binom{n-1}{k-1}=\binom{n-2}{k-1}+\binom{n-2}{k-2}=\binom{n-2}{k-1}
+\binom{n-3}{k-3}+ \binom{n-3}{k-2}$$
to rewrite the last inequality
in the following form
\begin{equation}
\binom{n-2}{k-1} \le \binom{n-2}{k}\,.
\end{equation}
An immediate verification shows that
(2.6) does not hold for the case when $n-k \le k$.
This implies that one of the conditions
(2.5) fails; therefore $n_{j}\le 1$.
Lemma 2.2 guarantees that $n_{j}>0$ and
we obtain the required.
$\square$
\begin{lemma}
The condition $n_{i}=2$
implies the
existence of unique line $L_{j}$ {\rm(}$j\ne
i${\rm)} contained in the plane $S_{i}$ and such that
$n_{j}=1$.
\end{lemma}
{\bf Proof.}  For the case when $n_{i}=2$
there exists
unique line
$L_{j}$ ($j\ne i$) contained
in $S_{i}$.
The trivial inclusion ${\cal R}'_{i} \subset {\cal R}'_{j}$
shows that $S_{j}\subset S_{i}$
and $0<n_{j}\le 2$.
If $n_{j}\ne 1$
then the planes $S_{i}$ and $S_{j}$ are
coincident.
Consider the plane
$S\in {\cal R}^{n-2}$ which does not contain
the lines $L_{i}$ and
$L_{j}$.
If a plane $U\in {\cal R}'$ does not belong to
${\cal R}'_{i}$ then the lines $L_{i}$ and $L_{j}$
are not contained in $U$
and
$U \in {\cal R}(S)$.
In other words, for each plane $U\in {\cal R}'$
the following two
cases can be realized:
\begin{enumerate}
\item[---]
$U \in {\cal R}'_{i}={\cal R}'_{j}\subset {\cal R}(S_{i})$,
\item[---]
$U \in {\cal R}(S)$.
\end{enumerate}
This implies the inclusion
${\cal R}' \subset
{\cal R}(S)\cup {\cal R}'_{i}$
showing that
\begin{equation*}
\begin{split}
\binom{n-1}{k-1}=\binom{n-2}{k-1}+\binom{n-2}{k-2}&\le
|{\cal R}'|\le\\
&|{\cal R}(S)|+
|{\cal R}(S_{i})|= \binom{n-2}{k}+\binom{n-2}{k-2}\,.
\end{split}
\end{equation*}
We obtain the inequality (2.6) but
for the case when $n-k \le k$ it does not hold
(see the proof of Lemma 2.4).
Thus
the equality $n_{j}=2$ fails and
we get $n_{j}=1$.
$\square$

Lemmas 2.2 and 2.4 show that if
$n_{i}\ne 1,2$ for some number $i$ then
$n({\cal R}')=n-1$.  Let us consider the case when $0< n_{i} \le 2$
for each $i=1,...,n$.
\begin{lemma}
If $0< n_{i} \le 2$ for each $i=1,...,n$
then
$$n({\cal R}')>k\; \mbox{ or }\; n({\cal R}')=1\,.$$
Moreover, the equality $n({\cal R}')=1$
holds if and only if there exists
a number $j$ such that ${\cal R}'={\cal R}(L_{j})$.
\end{lemma}
{\bf Proof.}
First of all note that for the case when $k=n-1$
our statement is trivial;
$\binom{n-1}{k-1}=n-1$ and ${\cal R}'$ is an
$R$-subset of ${\mathbb G}_{n-1}(V)$ containing not less than
$n-1$ elements.  If ${\cal R}'$ has $n$ elements then it is a
maximal $R$-set and $n({\cal R}')=n$.  If the number of elements
is equal to $n-1$ then there exists a plane $U\in {\cal R}$ such
that ${\cal R}'={\cal R}\setminus \{U\}$.  Consider unique line
$L_{j}$ which is not contained in $U$, it is trivial that ${\cal
R}'={\cal R}(L_{j})$.

Let $n-k\le k <n-1$. Fix a number $i$ such that $n_{i}=2$
and consider the plane
$S \in{\cal R}^{n-1}$ which does not contain
$L_{i}$.
It is not difficult
to see that the sets ${\cal R}(S)$ and ${\cal R}(S_{i})$ are
disjoint and
$$|{\cal R}'|=
|{\cal R}'\cap {\cal R}(S)|+|{\cal R}'_{i}|$$
(see the proof of Lemma 2.3). Thus
$$|{\cal R}'\cap {\cal R}(S)| =
|{\cal R}'| - |{\cal R}'_{i}|\ge
|{\cal R}'| - |{\cal R}(S_{i})|\ge$$
$$\binom{n-1}{k-1}-\binom{n-2}{k-2}=\binom{n-2}{k-1}\,.$$
In other words,
$${\cal R}'':={\cal R}'\cap {\cal R}(S)$$
is an $R$-subset of
${\mathbb G}_{k}(S)$ containing not less than
$\binom{(n-1)-1}{k-1}$ elements
(here ${\mathbb G}_{k}(S)$ is the Grassmannian of all
$k$-dimensional planes of the $(n-1)$-dimen\-sio\-n\-al vector
space $S$).

Note that if the set ${\cal R}''$ contains $\binom{n-2}{k-1}$
elements then
$$
|{\cal R}'_{i}|=|{\cal R}'|-|{\cal R}''|\ge
\binom{n-1}{k-1}-\binom{n-2}{k-1}=\binom{n-2}{k-2}\,.$$
Recall that the set ${\cal R}(S_{i})$
contains $\binom{n-2}{k-2}$ elements and
${\cal R}'_{i} \subset {\cal R}(S_{i})$.
This implies that for this case the set ${\cal R}'_{i}$
coincides with ${\cal R}(S_{i})$.
This remark will be exploited in what follows.

Now suppose that Lemma 2.6 holds for
$n-k<m$ (here $m$ is a natural number such that
$1<m<n-1$)
and consider the case when $n-k=m$.
By the inductive hypothesis and the remark made before Lemma
2.6 we have
$$n({\cal R}'')>k\; \mbox{ or }\; n({\cal R}'')=1\,.$$
The trivial inequality
$n({\cal R}')\ge n({\cal R}'')$
implies the fulfilment of the required statement
for the first case.

For the second case the inductive hypothesis
implies
the existence of a
number $j_{1}\ne i$ such that
\begin{equation} {\cal R}''={\cal
R}(S)\cap {\cal R}(L_{j_{1}})\,.
\end{equation}
Then the set
${\cal R}''$ contains $\binom{n-2}{k-1}$ elements and
${\cal R}'_{i}$ coincides with ${\cal R}(S_{i})$.
By Lemma 2.5 there exists
unique line $L_{j_{2}}$ ($j_{2}\ne
i$) contained in the plane $S_{i}$ and such that $n_{j_{2}}=1$.
It is easy to see that
$${\cal R}(L_{j_{2}})=
{\cal R}(S_{i})\cup ({\cal R}(S)\cap {\cal R}(L_{j_{2}}))\,.$$
Then the equation (2.7) and the equality
${\cal R}(S_{i})={\cal R}'_{i}$
show that
$${\cal R}(L_{j_{2}})={\cal R}'_{i}\cup{\cal R}''={\cal R}'$$
if $j_{1}=j_{2}$.

Consider the case when
$j_{1}\ne j_{2}$ and prove that
$n_{j}=1$ for each $j\ne i$.
For $j=j_{1}$ or $j_{2}$ it is trivial.
For $j\ne j_{1},j_{2}$ denote by  $S'_{j}$
the intersection
of all planes belonging to ${\cal R}''$ and
containing the line $L_{j}$. It is the two-dimensional plane
containing the lines $L_{j}$ and $L_{j_{1}}$.
The condition $j_{1}\ne j_{2}$ implies the existence of
a plane $U\in {\cal R}(S_{i})$
which contains $L_{j}$ and does not contain $L_{j_{1}}$.
The intersection of $U$
with $S'_{j}$ is the line $L_{j}$.
The equality ${\cal R}(S_{i})={\cal R}'_{i}$ guarantees that
the plane $U$ belongs to ${\cal R}'$.
Thus $n_{j}=1$ for each $j\ne i$ and $n({\cal R}')=n-1$.
$\square$

{\bf Proof of Theorem 2.1.}
It was noted above that we can restrict ourself
only to the case when $n-k \le k <n-1$.
Lemmas 2.2, 2.4 and 2.6
show that
we have to consider the following four cases.
\begin{enumerate}
\item[(i)] The inequality $n_{i}>0$ holds for any $i=1,...,n$
and there exists a number $j$ such that $n_{j} \ge 3$. Then
$n({\cal R}')=n-1$ (Lemma 2.4).
\item[(ii)] $0< n_{i} \le 2$ for each $i=1,...,n$ and
$n({\cal R}')>k$.
\item[(iii)] $0< n_{i} \le 2$ for each $i=1,...,n$ and
$n({\cal R}')=1$. Then ${\cal R}'$ coincides with some
${\cal R}(L_{j})$.
\item[(iv)] The condition
$n_{i}=0$ holds for some number $i$.  Then $n=2k$ and ${\cal
R}'={\cal R}(S)$, where $S$ is a plane belonging to the set
${\cal R}^{n-1}$ (Lemma 2.2).
\end{enumerate}

{\it Case }(i).
Lemma 2.3 states that
$n_{j} \le n-k$
and there exist $k-1$ numbers $i_{1},...,i_{k-1}$ such that
the plane $S_{j}$ does not contain the lines
$L_{i_{1}},...,L_{i_{k-1}}$. Denote by  $U$ the
$k$-dimensional plane
containing the lines
$L_{i_{1}},...,L_{i_{k-1}}$ and $L_{j}$.
For a number $i\ne j$ we have $n_{i}=1$ and
the intersection $U \cap S_{j}$ coincides
with the line $L_{j}$.
Therefore
the $R$-set ${\cal R}' \cup \{U\}$ is exact
and $\deg({\cal R}')=1$.

{\it Case }(ii).
Consider all numbers
$i_{1},...,i_{m}$ such that
$$n_{i_{1}}=...=n_{i_{m}}=2\;.$$
Then $n({\cal R}')=n-m$ and the condition $n({\cal R}')>k$ shows that
$m <n-k \le k$.
Denote by $S$ the plane containing
the planes $S_{i_{1}},...,S_{i_{m}}$ and having
the smallest dimension.
The dimension of
$S$ is not greater than $2m$
and
$$n-2m =(n-m)-m> k-m> 0\,.$$
This implies the existence of
$k-m$ numbers $j_{1},...,j_{k-m}$ such that
$$n_{j_{1}}=...=n_{j_{k-m}}=1$$
and $S$  does not contain  the lines
$L_{j_{1}},...,L_{j_{k-m}}$.
Denote by $U$ the
$k$-dimensio\-nal plane
containing the lines
$$L_{i_{1}},...,L_{i_{m}}, L_{j_{1}},...,L_{j_{k-m}}\,.$$
For any number $p=1,...,m$ the plane
$S_{i_{p}}$ is generated by two lines, one of them coincides with
$L_{i_{p}}$, other line $L_{q(p)}$ satisfies the
following conditions $n_{q(p)}=1$ (Lemma 2.5) and
$$q(p)\ne j_{1},...,j_{k-m}\,.$$ Then the intersection of the
plane $U$ with each $S_{i_{p}}$ is the line $L_{i_{p}}$ and the
$R$-set ${\cal R}' \cup \{U\}$ is exact.  We obtain the equality
$\deg({\cal R}')=1$.

{\it Cases }(iii)  and (iv)
were considered in Examples 2.1 and 2.2, respectively.
$\square$

\subsection{Proof of Theorem 1.1 for the general case}

First of all note that
the case $n-k<k$ can be reduced
to the case $n-k>k$.
It is not difficult to see that the following statements
hold true.
\begin{enumerate}
\item[---] For any two transformations
$f$ and $g$ of ${\mathbb G}_{k}(V)$ and
${\mathbb G}_{n-k}(V)$
induced by collineations and each
non-degenerate sesquilinear form $\Omega$
the bijections $f_{k\,n-k}(\Omega)f$
and $gf_{k\,n-k}(\Omega)$
are defined by sesquilinear forms.
\item[---] For any  bijection
defined by a sesquilinear form
the inverse bijection is
defined by a sesquilinear form.
\item[---] The composition of two bijections defined
by sesquilinear forms is a transformation
induced by some collineation.
\end{enumerate}
For each regular transformation
$f$ of ${\mathbb G}_{k}(V)$ and a non-degenerate sesquilinear form
$\Omega$ on $V$
$$g=f_{k\,n-k}(\Omega)fF_{n-k\,k}(\Omega)$$
is a regular transformation of ${\mathbb G}_{n-k}(V)$.
The statements given above show that
$f$ is induced by a collineation if and only if
$g$ is induced by a collineation.

The similar arguments can
be used to prove Corollary 1.1, since for each regular bijection
$f$ of ${\mathbb G}_{k}(V)$ onto ${\mathbb G}_{n-k}(V)$ and a
non-degenerate sesquilinear form $\Omega$ on $V$ the composition
$F_{n-k\,k}(\Omega)f$ is a regular
transformation of ${\mathbb G}_{k}(V)$.

Let $f$ be a regular transformation of ${\mathbb G}_{k}(V)$ and
$1<k\le n-k$. Let also $U$ and $U'$ be $k$-dimensional adjacent
planes in $V$. We want to show  that the planes $f(U)$ and $f(U')$
are adjacent. Consider a maximal $R$-set ${\cal R}$ containing $U$
and $U'$.  Then the maximal $R$-set $f({\cal R})$ contains $f(U)$
and $f(U')$; in what follows this set will be denoted by ${\cal
R}_{f}$.  For planes $S$ and $S'$ belonging to  ${\cal R}^{m}$ and
${\cal R}^{m}_{f}$ (here ${\cal R}^{m}_{f}$ is the maximal
$R$-subset of ${\mathbb G}_{m}(V)$ associated with ${\cal R}_{f}$)
denote by ${\cal R}(S)$ and ${\cal R}_{f}(S')$ the sets of all
planes $U\in {\cal R}$ and $U'\in {\cal R}_{f}$ incident to $S$
and $S'$, respectively.

\begin{lemma}
For any plane $S\in{\cal R}^{n-1}$ the
following statements are fulfilled:
\begin{enumerate}
\item [{\rm (i)}] if $n \ne 2k$ then there exists
a plane $S'\in {\cal R}^{n-1}_{f}$
such that
\begin{equation}
f({\cal R}(S))={\cal R}_{f}(S')\,;
\end{equation}
\item
[{\rm (ii)}] if $n=2k$ then there exists
a plane $S'\in  {\cal R}^{m}_{f}$
such that
$m=1$ or $n-1$ and the equality
{\rm(2.8)} holds true.
\end{enumerate}
\end{lemma}

{\bf Proof.}
It is a direct consequence of Theorem 2.1 and
Lemma 2.1.
$\square$

\begin{lemma}
For any plane $S\in{\cal R}^{k+1}$
the following statements are fulfilled:
\begin{enumerate}
\item [{\rm (i)}] if $n \ne 2k$ then there exists
a plane $S'\in {\cal R}^{k+1}_{f}$
such that the equality
{\rm(2.8)} holds true;
\item[{\rm (ii)}] if $n=2k$ then there exists
a plane $S'\in {\cal R}^{m}_{f}$
such that
$m=k-1$ or $k+1$ and the equality
{\rm(2.8)} holds true.
\end{enumerate}
\end{lemma}

{\bf Proof.}
Let $L_{1},...,L_{n}$ and $L'_{1},...,L'_{n}$
be planes generating the sets
${\cal R}^{1}$ and
${\cal R}^{1}_{f}$, respectively.
Denote by $S_{i}$ and $S'_{i}$
the planes belonging to
${\cal R}^{n-1}$ and
${\cal R}^{n-1}_{f}$ and such that
$S_{i}$  does not contain $L_{i}$ and $S'_{i}$
does not contain
$L'_{i}$.
Lemma 2.7 states that if $n\ne 2k$ then for each
$i=1,...,n$ there exists a number $j_{i}$ such that
\begin{equation}
f({\cal R}(S_{i}))={\cal R}_{f}(S'_{j_{i}})\,.
\end{equation}
Let us prove the statement (i) for the case when
$S$ is the $(k+1)$-dimensional plane containing
the lines $L_{1},...,L_{k+1}$
(for other planes belonging to
the set ${\cal R}^{k+1}$ the proof is similar).
We have $S=\cap^{n}_{i=k+2}S_{i}$ and
$${\cal R}(S)=\bigcap^{n}_{i=k+2}{\cal R}(S_{i})\,.$$
The last equality and the equation (2.9) show that
$$f({\cal R}(S))=
\bigcap^{n}_{i=k+2}{\cal R}_{f}(S'_{j_{i}})
={\cal R}_{f}(S')\,;$$
where $S'$ is the $(k+1)$-dimensional plane
containing the lines $L'_{j_{1}},...,L'_{j_{k+1}}$.

If $n=2k$ then
by Lemma 2.7 the following two cases can be realized:
\begin{enumerate}
\item [(a)] there exists a number $j_{1}$
such that the equality (2.9)
holds for $i=1$;
\item [(b)] there
exists a number $j$ such that
$f({\cal R}(S_{1}))={\cal R}_{f}(L'_{j})$.
\end{enumerate}
First of all show that for the case (a)
there
exist numbers $j_{2},...,j_{n}$ such that  the equality (2.9)
holds for any $i=2,...,n$ .  Then the proof of the statement (ii)
will be similar to the proof of the statement (i) given above.

Suppose that
$$f({\cal R}(S_{i}))={\cal R}_{f}(L'_{j_{i}})$$
for some numbers $i$ and $j_{i}$ and
consider the $(n-2)$-dimensional plane
${\hat S} = S_{1} \cap S_{i}$.  Then
$${\cal R}({\hat S})=
{\cal R}(S_{1}) \cap {\cal R}(S_{i})$$
and
$$f({\cal R}({\hat S}))=
{\cal R}_{f}(S'_{j_{1}})\cap {\cal R}_{f}(L'_{j_{i}})\,.$$
The first set contains $\binom{2k-2}{k}$ elements,
the second set contains
$\binom{2k-2}{k-1}$ elements.
However
$$\binom{2k-2}{k} \ne \binom{2k-2}{k-1}\,.$$
Thus our hypothesis
fails and the statement (ii) is proved for the case (a).

Show that the case (b) can be reduced to the case (a).
Consider a
non-degenerate sesquilinear form $\Omega$ on $V$
and the regular transformation
$g=F_{k\,k}(\Omega)$.
It is easy to see that the regular
transformation $gf$ satisfies the conditions of the case (a).
Then the maximal $R$-subset of ${\mathbb G}_{k+1}(V)$
associated with ${\cal R}_{gf}$
contains a plane
$S''$
such that
$$gf({\cal R}(S))= {\cal R}_{gf}(S'')\,.$$
Then the plane
$$S'=(F_{k-1\,k+1}(\Omega))^{-1}(S'')\in {\cal R}^{k-1}_{f}$$
satisfies the required condition.
$\square$

Now we can prove Theorem 1.1. The planes $U$ and $U'$ are adjacent
and there exists a $(k+1)$-dimensiomal plane $S \in {\cal
R}^{k+1}$ containing them. By Lemma 2.8 the planes $f(U)$ and
$f(U')$ are adjacent. The inverse transformation $f^{-1}$ is
regular and the similar arguments show that the planes $f^{-1}(U)$
and $f^{-1}(U')$ are adjacent too. We have proved that $f$
preserves the distance between planes and the required statement
is a direct consequence of the Chow Theorem.

\subsection{Remark on exact $R$-sets}

Theorem 2.1 states that each $R$-set ${\cal R}\subset {\mathbb
G}_{k}(V)$ containing not less than $\binom{n-1}{k-1}+1$ elements
if $n-k\le k$ and $\binom{n-1}{k}+1$ elements if $n-k\ge k$
satisfies the condition $\deg {\cal R}\le 1$. It is natural to ask
what is the minimal number $s^{n}_{k}$ such that each $R$-subset
of ${\mathbb G}_{k}(V)$ containing more than $s^{n}_{k}$ elements
is exact? The following statement gives the answer on the
question.
\begin{theorem}
If $1<k<n-1$ then each $R$-subset
of ${\mathbb G}_{k}(V)$ containing more
than
$$s^{n}_{k}:=\binom{n-1}{k}+\binom{n-2}{k-2}$$
elements is exact.
Moreover, there exists a non-exact $R$-subset
of ${\mathbb G}_{k}(V)$ containing
$s^{n}_{k}$ elements.
\end{theorem}
Theorem 2.2 is
not connected with Theorem 1.1, but its proof is not
complicated and we give it here.
It will be based on the terms introduced
in Subsection 2.1.
\begin{exmp}{\rm
Let ${\cal R}$ be a maximal
$R$-subsets of ${\mathbb G}_{k}(V)$.
Consider two planes
$S \in {\cal R}^{n-1}$ and $S'\in {\cal R}^{2}$
such that $S$ does not contain $S'$.
Suppose that
$${\cal R}'={\cal R}(S)\cup {\cal R}(S')\,.$$
The sets ${\cal R}(S)$ and ${\cal R}(S')$ are disjoint and
$$|{\cal R}'|=|{\cal R}(S)|+|{\cal R}(S')|=
\binom{n-1}{k}+\binom{n-2}{k-2}=s^{n}_{k}\,.$$
There exists unique line $L_{i}$ which is not contained in
the plane $S$
(this line is contained in $S'$).  Then $n_{j}=1$ if $j\ne i$
and $S_{i}=S'$. Thus $n_{i}=2$
and the $R$-set ${\cal R}'$ is not exact.
Consider a
$k$-dimensional plane
$U\in {\cal R}$ which contains the line $L_{i}$ and does not
contain the plane $S'=S_{i}$.  The intersection $U\cap S_{i}$
coincides with $L_{i}$. This implies that the $R$-set ${\cal
R}'\cup \{U\}$ is exact and $\deg({\cal R}')=1$.
}\end{exmp}
Theorem 2.2 is a direct consequence of the following statement.
\begin{theorem}
If an $R$-set ${\cal R}'\subset {\mathbb G}_{k}(V)$, $1<k<n-1$,
contains not less than $s^{n}_{k}$ elements then
$\deg({\cal R}') \le 1$ and the equality
$\deg({\cal R}')=1$ holds if and only if
${\cal R}'$ is the set considered in Example
{\rm 2.3}.
\end{theorem}
{\bf Proof.}
If the $R$-set ${\cal R}'$ is not exact then there exists a
number $i$ such that $n_{i}\ne 1$.  Consider unique
plane $S\in {\cal R}^{n-1}$
which does not contain
the line $L_{i}$. It was noted above
that
${\cal R}'$ can be represented as
the union of the two disjoint sets
${\cal R}(S)\cap {\cal R}'$ and ${\cal R}'_{i}$.
We have
$$ |{\cal R}(S)\cap {\cal R}'|\le |{\cal R}(S)|=
\binom{n-1}{k}$$
and
$$|{\cal R}'_{i}=R(S_{i})\cap {\cal R}'|\le
|R(S_{i})|=\binom{n-n_{i}}{k-n_{i}}\le \binom{n-2}{k-2}$$
(see Remark 2.1).
Then the inequality $$|{\cal R}'|=|{\cal
R}(S)\cap {\cal R}'|+|{\cal R}'_{i}|
\le \binom{n-1}{k}+\binom{n-2}{k-2}
=s^{n}_{k}$$ shows that the condition
$|{\cal R}'| \ge s^{n}_{k}$
holds if and only if
$${\cal R}(S)\cap {\cal R}'= {\cal R}(S)\,,\;\;
{\cal R}'_{i}={\cal R}(S_{i})$$ and
$n_{i}=2$.  $\square$

\section{Transformations of ${\mathfrak I}_{k\,n-k}(V)$
and automorphisms of classical groups}

\setcounter{theorem}{0}
\setcounter{lemma}{0}
\setcounter{prop}{0}
\setcounter{exmp}{0}
\setcounter{rem}{0}
\setcounter{equation}{0}

\subsection{Transformations of ${\mathfrak I}_{k\,n-k}(V)$
preserving the commutativity and the adjacency}

Denote by ${\mathbb G}_{k\,n-k}(V)$
the set of all pairs
$$(U, S)\in {\mathbb G}_{k}(V)\times{\mathbb
G}_{n-k}(V)$$
such that
$$U+S=V\,.$$
We say that ${\cal R}\subset {\mathbb G}_{k\,n-k}(V)$
is an $R$-{\it set} if there exists
a base for $V$ such that for each
pair $(U, S)$ belonging to ${\cal R}$
the planes $U$ and $S$ contain $k$ and
$n-k$ vectors from this base;
i.e. there exists a coordinate system
such that for all $(U, S)\in{\cal R}$
the planes $U$ and $S$ are
coordinate planes for this system.
Any base and any coordinate system satisfying
that condition will be called {\it associated}
with ${\cal R}$.

An $R$-set ${\cal R}$ is called
{\it maximal} if any $R$-set containing ${\cal R}$
coincides with it. It is trivial that
an $R$-subset of ${\mathbb G}_{k\,n-k}(V)$
is maximal if and only if it contains
$\binom{n}{k}$ elements.

We say that a transformation of ${\mathbb G}_{k\,n-k}(V)$
is {\it regular} if it preserves the
class of $R$-sets.
A bijection $f$ of
${\mathbb G}_{k\,n-k}(V)$ onto
${\mathbb G}_{n-k\,k}(V)$
is called {\it regular} if $f$ and $f^{-1}$
map each $R$-set to an $R$-set.
For the case when $m\ne k,n-k$
there are not regular bijections of
${\mathbb G}_{k\,n-k}(V)$ onto
${\mathbb G}_{m\,n-m}(V)$.

If $F$ is a field with characteristic other than two
then denote by $\pi$
the bijection
$$\sigma \to (U_{+}(\sigma),U_{-}(\sigma))$$
of ${\mathfrak I}_{k\,n-k}(V)$ onto
${\mathbb G}_{k\,n-k}(V)$.

\begin{prop}
For a set ${\cal R}\subset {\mathbb G}_{k\,n-k}(V)$
the following two conditions are equivalent:
\begin{enumerate}
\item[{\rm (i)}] ${\cal R}$ is an $R$-set,
\item[{\rm (ii)}]involutions belonging to the set $\pi^{-1}({\cal R})$
commute.
\end{enumerate}
\end{prop}

The proof of Proposition 3.1 will be based on the following
well-known lemma.

\begin{lemma}
An involution $s$
and a linear transformation $f$
commute if and only if $f$ preserves
$U_{+}(s)$ and $U_{-}(s)$.
\end{lemma}

{\bf Proof.}
The implication
${\rm (i)}\Rightarrow {\rm (ii)}$
is a simple consequence of Lemma 3.1

Suppose that ${\cal R}$ is a subset of ${\mathbb G}_{k\,n-k}(V)$
such that any two elements of $\pi^{-1}({\cal R})$ commute.
If the number of
elements $|{\cal R}|$ is finite
then for each set ${\cal R}'\subset {\cal R}$
define
$$U_{{\cal R}'}=\left(
\bigcap_{s\in \pi^{-1}({\cal R}')}U_{+}(s) \right) \bigcap
\left(
\bigcap_{s\in \pi^{-1}({\cal R}\setminus{\cal R}')}U_{-}(s)
\right).$$
Some of these intersections are non-zero, we denote them
by $U_{1},...,U_{i}$. Then
$$U_{1}+...+U_{i}=V\,,$$
$U_{i}\cap U_{j}=\{0\}$ if $i\ne j$ and there exists
a coordinate system such that $U_{1},...,U_{i}$
are coordinate planes for this system.
This statement can be obtained as a consequence of
Lemma 3.1 by the induction on
$|{\cal R}|$.
It is trivial that for any
$(U,S)\in {\cal R}$
the planes $U$ and $S$ are coordinate planes for our system
and ${\cal R}$ is an $R$-set.

For the general case the arguments given above
show that each finite subset of ${\cal R}$
contains not greater than $\binom{n}{k}$
elements. Thus the set ${\cal R}$ is finite and
we get the required.
$\square$

Proposition 3.1 shows that a transformation
$f$ of ${\mathbb G}_{k\,n-k}(V)$ is
regular if and only if $\pi^{-1} f \pi$
is a transformation of ${\mathfrak I}_{k\,n-k}(V)$
preserving the commutativity.
\begin{exmp}{\rm
For a collineation $g$ of $V$
denote by $g_{k}$ and $g_{n-k}$
the transformations of ${\mathbb G}_{k}(V)$
and ${\mathbb G}_{n-k}(V)$ induced by it.
Then
$$(U,S)\to (g_{k}(U),g_{n-k}(S))\,,$$
is a regular transformation of
${\mathbb G}_{k\,n-k}(V)$.
The respective transformation
of ${\mathfrak I}_{k\,n-k}(V)$ maps
each involution $\sigma$
to $g\sigma g^{-1}$.
}\end{exmp}
\begin{exmp}{\rm
Recall that any non-degenerate sesquilinear
form $\Omega$ defines the regular bijection
$f_{k\,n-k}(\Omega)$ of
${\mathbb G}_{k}(V)$ onto
${\mathbb G}_{n-k}(V)$.
The map
$$(U,S)\to (f_{n-k\,k}(\Omega)(S),f_{k\,n-k}(\Omega)(U))$$
is a regular transformation of ${\mathbb G}_{k\,n-k}(V)$;
denote it by $f_{k\,n-k}(\Omega)$.
It is not difficult to see that
$$\pi^{-1}f_{k\,n-k}(\Omega)\pi$$
transfers each $\sigma$ to
$h^{-1}\check{\sigma}h$, where
$h$ is the correlation defined by our form
$\Omega$.
}\end{exmp}
\begin{exmp}{\rm
Consider the bijection $i_{k\, n-k}$
of ${\mathbb G}_{k\,n-k}(V)$
onto ${\mathbb G}_{n-k\,k}(V)$
which
transfers $(U,S)$
to $(S,U)$. Then $i_{k\, n-k}$ is regular
and $\pi^{-1} i_{k\, n-k} \pi$
coincides with the bijection of ${\mathfrak I}_{k\,n-k}(V)$
onto ${\mathfrak I}_{n-k\,k}(V)$
transferring $\sigma$ to $-\sigma$.
}\end{exmp}
We say that two pairs $(U,S)$ and $(U',S')$
belonging to ${\mathbb G}_{k\, n-k}(V)$
are {\it adjacent} if
one of the following conditions holds true:
\begin{enumerate}
\item[---] $U=U'$ and the planes $S$ and $S'$
are adjacent,
\item[---] $S=S'$ and the planes $U$ and $U'$
are adjacent.
\end{enumerate}
Two involutions $\sigma$ and $\sigma'$
will be called {\it adjacent} if the pairs
$\pi(\sigma)$ and $\pi(\sigma')$
are adjacent.
It is trivial that the transformations
considered in Examples 3.1 -- 3.3
preserve the adjacency.

An immediate verification shows that
two involutions $\sigma$ and $\sigma'$
are adjacent if and only if their composition
$\sigma\sigma'$ is a transvection
(recall that a linear transformation
$g\in {\mathfrak S}{\mathfrak L}(V)$ is called a transvection if
the dimension of $\ker(Id -g)$ is equal to $n-1$).
\begin{theorem}
If $n\ge 3$ and the characteristic of the field $F$
is not equal to $2$ then the
following two statements are fulfilled:
\begin{enumerate}
\item[---]  for the case when $n\ne 2k$
each transformation of
${\mathfrak I}_{k\,n-k}(V)$ preserving the commutativity
and the adjacency
is defined by a collineation $($Example {\rm 3.1}$)$
or a correlation $($Example {\rm 3.2}$)$;
\item[---] if $n=2k$ then
for any transformation $f$ of
${\mathfrak I}_{k\,n-k}(V)$ preserving the commutativity
and the adjacency
one of the transformations $f$ or
$-f$ $($Example {\rm 3.3}$)$
is defined by a collineation or a correlation.
\end{enumerate}
\end{theorem}
\begin{cor}
If $n\ge 3$ and the characteristic of the field $F$
is not equal to $2$ then the
following two statements hold true:
\begin{enumerate}
\item[---]  for the case when $n\ne 2k$
each transformation of
${\mathfrak I}_{k\,n-k}(V)$ preserving the commutativity
the adjacency
can be extended to
an automorphism of the group generated by
$(k,n-k)$-involutions;
\item[---] if $n=2k$ then
for any transformation $f$ of
${\mathfrak I}_{k\,n-k}(V)$ preserving the commutativity
and the adjacency
one of the transformations $f$ or
$-f$ can be extended to
an automorphism of the group generated by
$(k,n-k)$-involutions.
\end{enumerate}
\end{cor}
For $k=1,n-1$ each transformation of ${\mathfrak I}_{k\,n-k}(V)$
preserving the commutativity preserves the adjacency. For the case when
$n=2k$ this statement fails (Example 3.4). For other cases it  is not
proved yet; we have not an analogy of Mackey's lemma (see [M] or [D2]).

Proposition 3.1 shows that Theorem 3.1 can be reformulated in the
following form.
\begin{theorem}
If $n\ge 3$ then the
following two statements hold true:
\begin{enumerate}
\item[---]  for the case when $n\ne 2k$
each regular transformation of
${\mathbb G}_{k\,n-k}(V)$ preserving the adjacency
is
induced by a collineation
or defined by a non-degenerate
sesquilinear form;
\item[---] if $n=2k$ then
for any regular transformation $f$ of
${\mathbb G}_{k\,n-k}(V)$ preserving the adjacency
one of the transformations $f$ or
$i_{k\,k}f$ is
induced by a collineation
or defined by a non-degenerate
sesquilinear form.
\end{enumerate}
\end{theorem}
\begin{exmp}{\rm
Suppose that $n=2k$ and consider a set
${\cal G}\subset {\mathbb G}_{k\,k}(V)$
satisfying the following condition:
for each $(U,S)\in {\cal G}$ the pair
$(S,U)$ belongs to ${\cal G}$.
Define
$$i_{{\cal G}}(U,S)=(S,U)\;\mbox{ if }\;(U,S)\in {\cal G}$$
and
$$i_{{\cal G}}(U,S)=(U,S)\;\mbox{ if }\;(U,S)\notin {\cal G}\,.$$
Each maximal $R$-subset of ${\mathbb G}_{k\,k}(V)$
containing $(U,S)$ contains $(S,U)$.
This implies that the transformation $i_{{\cal G}}$ is regular
but for the case when ${\cal G}\ne \emptyset$ and
${\mathbb G}_{k\,k}(V)$ it does not preserve the adjacency.
}\end{exmp}

\subsection{Proof}

{\bf First step.}
For a $k$-dimensional plane $U$ and
an $(n-k)$-dimensional plane $S$ denote by ${\cal X}(U)$ and
${\cal X}(S)$ the sets of all pairs $(U,\cdot)$ and $(\cdot, S)$
belonging to ${\mathbb G}_{k\,n-k}(V)$. It is trivial
that
the intersection
${\cal X}(U)\cap {\cal X}(S)$ is not empty
if and only if
$U+S=V$; for this case it consists only of the pair
$(U,S)$.
For other $k$-dimensional plane $U'$
we have
${\cal X}(U)\cap {\cal X}(U')=\emptyset$.

Now for arbitrary plane $P\subset V$ denote by
${\cal Y}(U, P)$ the set of all pairs
$(U,S')\in {\cal X}(U)$ such that $S'$ is incident
to $P$. Denote also by ${\cal Y}(S, P)$
the set of all pairs
$(U',S)\in {\cal X}(S)$ such that $U'$ is incident
to $P$. Note that for some cases the sets
${\cal Y}(U, P)$ and ${\cal Y}(U, P)$ may be empty
(for example, if $P$ is a plane contained in $U$ or $S$).

If $P$ is an $(n-k\pm 1)$-dimensional
plane then  any two pairs belonging to
${\cal Y}(U, P)$ are adjacent. The set
${\cal Y}(S, P)$ satisfies the similar condition
if the dimension of $P$ is equal to $k\pm 1$.
In what follows we will use the following statement:
if ${\cal G}$ is a subset of ${\mathbb G}_{k\, n-k}(V)$
any two elements of which are adjacent then
there exists a $k$-dimensional plane $U$
or an $(n-k)$-dimensional plane $S$ such that
${\cal G}$ is contained in
${\cal X}(U)$ or ${\cal X}(S)$ (the proof is trivial).

{\bf Second step.}
Let $f$ be a transformation of
${\mathbb G}_{k\,n-k}(V)$ preserving the adjacency.
Show that {\it for
any $k$-dimensional plane $U$ there exists
a $k$-dimensional plane $U_{f}$ or
an $(n-k)$-dimensional plane $S_{f}$ such that
the set $f({\cal X}(U))$ coincides with ${\cal X}(U_{f})$
or ${\cal X}(S_{f})$.}

Let us fix two adjacent pairs $(U,S)$ and $(U,S')$.
Then $P=S\cap S'$ is an $(n-k-1)$-dimensional plane
and any two pairs belonging to $f({\cal Y}(U,P))$
are adjacent.
Therefore there exists
a plane $T$ the dimension of which is equal
to $k$ or $n-k$ and such that
$$f({\cal Y}(U,P))\subset {\cal X}(T)\,.$$
Now consider a pair $(U, {\hat S})$ adjacent to
$(U,S)$. The planes $S$ and ${\hat S}$ are contained
in some $(n-k+1)$-dimensional plane $P'$.
There exists
a plane $T'$ such that $\dim T'=k$ or $n-k$
and
$$f({\cal Y}(U,P'))\subset {\cal X}(T')\,.$$
The inclusions $P\subset S\subset P'$ guarantees that the set
$${\cal Y}(U,P)\cap {\cal Y}(U,P')$$
contains not less than two
elements. Then
${\cal X}(T)\cap {\cal X}(T')$
satisfies the similar condition.
It means
that the planes $T$ and $T'$ are coincident
(see First step).

Thus  for any pair $(U, {\hat S})$
satisfying the condition $d(S,{\hat S})=1$
the pair $f(U, {\hat S})$ belongs to
the set ${\cal X}(T)$.
Suppose that the similar statement holds for
the case when $d(S,{\hat S})<i$, $i>1$ and consider
the case $d(S,{\hat S})=i$.

By the definition of the distance
there exists
a pair $(U,{\check S})\in {\cal X}(U)$
such that $d(S,{\check S})=i-1$ and
the planes ${\check S}$ and ${\hat S}$
are adjacent. Denote by $P''$
the $(n-k+i-1)$-dimensional
plane containing $S$ and ${\check S}$.
By our hypothesis
$$f({\cal Y}(U,P''))\subset {\cal X}(T)\,.$$
For the $(n-k-1)$-dimensional plane
${\hat P}={\hat S}\cap {\check S}$
there exists a plane ${\hat T}$
the dimension of which
is equal to $k$ or $n-k$ and such that
$$f({\cal Y}(U,{\hat P}))\subset {\cal X}({\hat T})\,.$$
Then
${\hat P}\subset {\check S}\subset P''$
and the set
$${\cal Y}(U,P'')\cap {\cal Y}(U,{\hat P})$$
contains not less than two elements. Thus ${\hat
T}$ coincides with $T$ and
we get the required.

{\bf Third step.}
For the transformation $f$ introduced above
consider the set ${\cal U}$
of all $k$-dimensional planes $U$
such that
\begin{equation}
{\cal X}(U)=f({\cal X}(U'))
\end{equation}
for some $k$-dimensional plane $U'$.
Consider also the set ${\cal S}$  of all
$(n-k)$-dimensional planes $S$ such that
\begin{equation}
{\cal X}(S)=f({\cal X}(U''))
\end{equation}
for some $k$-dimensional plane $U''$.

Let $U$ and $S$ be planes belonging to the sets
${\cal U}$ and ${\cal S}$. Let also $U'$ and
$U''$ be $k$-dimensional planes satisfying
the conditions (3.1) and (3.2).
Then ${\cal X}(U')\cap {\cal X}(U'')=\emptyset$.
Therefore ${\cal X}(U)\cap {\cal X}(S)=\emptyset$.
The last equality holds only for the case when
\begin{equation}
\dim U\cap S\ge 1\,.
\end{equation}
In other words,
the following statements are fulfilled:
\begin{enumerate}
\item[---]
if a $k$-dimensional plane $U$ belongs to ${\cal U}$ then the
inequality (3.3) holds for each $S\in {\cal S}$,
\item[---] if an
$(n-k)$-dimensional plane $S$ belongs to ${\cal S}$ then the
inequality (3.3) holds for each $U\in {\cal U}$.
\end{enumerate}
For arbitrary pair $(U,S)\in
{\mathbb G}_{k\,n-k}(V)$ one of the planes $U$ or $S$ belongs to
${\cal U}$ or ${\cal S}$, respectively.  Thus one of
these sets is empty.  If the transformation $f$ is regular
then one of the following cases is realized:
\begin{enumerate}
\item[---] for each $k$-dimensional plane $U$
there exists a $k$-dimensional plane $U_{f}$ such that
$$f({\cal X}(U))={\cal X}(U_{f})\,,$$
then $U\to U_{f}$ is a regular
transformation of ${\mathbb G}_{k}(V)$;
\item[---] for each
$k$-dimensional plane $U$ there exists an $(n-k)$-dimensional
plane $S_{f}$ such that
$$f({\cal X}(U))={\cal X}(S_{f})\,,$$
then $U\to S_{f}$ is a regular bijection of ${\mathbb G}_{k}(V)$
onto ${\mathbb G}_{n-k}(V)$.
\end{enumerate}
Theorem 1.1 and Corollary 1.1 give the required.

\end{document}